\documentclass[12pt,reqno]{amsart}
\usepackage{amsmath}
\usepackage{amsthm}

\theoremstyle{plain}
\newtheorem{theorem}{Theorem}

\newtheorem{corollary}[theorem]{Corollary}

\theoremstyle{definition}
\newtheorem{definition}[theorem]{Definition}

\theoremstyle{remark}
\newtheorem*{remark}{Remark}

\numberwithin{equation}{section}
\numberwithin{theorem}{section}
\numberwithin{figure}{section}

\DeclareMathOperator{\tr}{Tr}

\newcommand{\qbin}[2]{\genfrac{[}{]}{0pt}{}{#1}{#2}}
\newcommand{\sumtop}[2]{\genfrac{}{}{0pt}{2}{#1}{#2}}

\def\bK{{\mathbf K}} \def\bz{{\mathbf z}} 
\def\ba{{\mathbf a}} \def\bb{{\mathbf b}}
\def\bm{{\mathbf m}} \def\bn{{\mathbf n}} 
\def\bfe{{\mathbf e}}
\def\bu{{\mathbf u}} \def\bp{{\mathbf p}} 
\def\bt{{\mathbf t}} \def\bz{{\mathbf z}} 
\def\bM{{\mathbf M}} \def\bN{{\mathbf N}}
\def\bQ{{\mathbf Q}}

\newcommand{\wh}[1]{\widehat{#1}}
\newcommand{\eql}{=} 
\newcommand{\qn}[1]{(q)_{#1}} \newcommand{\qny}[2]{({#1}q)_{#2}}
\newcommand{\id}{{1\!\!1}}
\newcommand{\ZZ}{{\mathbb Z}}   \newcommand{\QQ}{{\mathbb Q}}  
\newcommand{\NN}{{\mathbb N}}   
\newcommand{\la}{\lambda}     \newcommand{\bla}{\boldsymbol{\la}}
\newcommand{\sln}{\mathfrak{sl}_{n+1}} 
\newcommand{\whsln}{\wh{\mathfrak{sl}}_{n+1}}
\newcommand{\fg}{\mathfrak{g}}
\newcommand{\bep}{\boldsymbol{\epsilon}}
\newcommand{\examtw}{\begin{picture}(6,4)(0,-4)
\thicklines
\put(0,0){\line(1,0){3}}
\put(0,-3){\line(1,0){3}}
\put(0,0){\line(0,-1){3}}
\put(3,0){\line(0,-1){3}}
\thinlines
\multiput(0,0)(0,-1){2}{\line(1,0){6}}
\multiput(0,-2)(0,-1){2}{\line(1,0){4}}
\put(0,-4){\line(1,0){2}}
\multiput(0,0)(1,0){3}{\line(0,-1){4}}
\multiput(3,0)(1,0){2}{\line(0,-1){3}}
\multiput(5,0)(1,0){2}{\line(0,-1){1}}\end{picture}} 

\newcommand{\mnab}{\begin{picture}(24,14)(0,0)
\put(6,12){\line(1,0){14}}
\put(6,6){\line(1,0){10}}
\put(6,12){\line(0,-1){12}}
\put(16,12){\line(0,-1){6}}
\put(6,0){\line(1,1){6}}
\put(16,8){\line(1,1){4}}

\put(7,12){\line(0,-1){3.5}}
\put(8,12){\line(0,-1){2.5}}
\put(9,12){\line(0,-1){1.5}}
\put(10,12){\line(0,-1){0.5}}

\put(6,11){\line(1,0){3.5}}
\put(6,10){\line(1,0){2.5}}
\put(6,9){\line(1,0){1.5}}
\put(6,8){\line(1,0){0.5}}

\put(9,13){\vector(-1,0){3}}
\put(13,13){\vector(1,0){3}}
\put(5,10){\vector(0,1){2}}
\put(5,8){\vector(0,-1){2}}

\put(8,6.5){\vector(-1,0){2}}
\put(10,6.5){\vector(1,0){2}}
\put(15.5,9){\vector(0,-1){1}}
\put(15.5,11){\vector(0,1){1}}

\put(11,13){\makebox(0,0){$m+a$}}
\put(5,9){\makebox(-0.5,0){$n+b$}}
\put(9,6.5){\makebox(0,0){$m$}}
\put(15.5,10){\makebox(0,0){$n$}}
\end{picture} }

\newcommand{\zezezeze}{\begin{picture}(1,1)(0,-1)
\put(0.5,-0.5){\circle*{0.2}} \end{picture}} 

\newcommand{\zeonzeon}{\begin{picture}(1.5,1)(0,-1)
\put(0,0){\line(1,0){1.5}}
\put(0,-1){\line(1,0){0.5}}
\put(0,0){\line(0,-1){1}}
\put(0.5,-1){\line(1,1){1}}  \end{picture}} 

\newcommand{\onononon}{\begin{picture}(2,2)(0,-2)
\put(0,0){\line(1,0){2}}
\put(0,-1){\line(1,0){1}}
\put(0,0){\line(0,-1){2}}
\put(1,0){\line(0,-1){1}}
\put(0,-2){\line(1,1){2}}  \end{picture}} 

\newcommand{\ontwontw}{\begin{picture}(3,3)(0,-3)
\put(0,0){\line(1,0){3}}
\multiput(0,-1)(0,-1){2}{\line(1,0){1}}
\put(0,0){\line(0,-1){3}}
\put(1,0){\line(0,-1){2}}
\put(0,-3){\line(1,1){3}}  \end{picture}} 

\newcommand{\onthonth}{\begin{picture}(4,4)(0,-4)
\put(0,0){\line(1,0){4}}
\multiput(0,-1)(0,-1){3}{\line(1,0){1}}
\put(0,0){\line(0,-1){4}}
\put(1,0){\line(0,-1){3}}
\put(0,-4){\line(1,1){4}}  \end{picture}} 

\newcommand{\zetwzetw}{\begin{picture}(2.5,2)(0,-2)
\put(0,0){\line(1,0){2.5}}
\put(0,-1){\line(1,0){0.2}}
\put(0,-2){\line(1,0){0.5}}
\put(0,0){\line(0,-1){2}}
\put(0.5,-2){\line(1,1){2}}  \end{picture}} 

\newcommand{\ononzeon}{\begin{picture}(2,1)(0,-1)
\put(0,0){\line(1,0){2}}
\put(0,-1){\line(1,0){1}}
\multiput(0,0)(1,0){2}{\line(0,-1){1}}
\put(1,-1){\line(1,1){1}}  \end{picture}} 

\newcommand{\ontwzetw}{\begin{picture}(3,2)(0,-2)
\put(0,0){\line(1,0){3}}
\multiput(0,-1)(0,-1){2}{\line(1,0){1}}
\multiput(0,0)(1,0){2}{\line(0,-1){2}}
\put(1,-2){\line(1,1){2}}  \end{picture}}

\newcommand{\twfotwfo}{\begin{picture}(6,6)(0,-6)
\put(0,0){\line(1,0){6}}
\multiput(0,-1)(0,-1){4}{\line(1,0){2}}
\put(0,0){\line(0,-1){6}}
\multiput(1,0)(1,0){2}{\line(0,-1){4}}
\put(0,-6){\line(1,1){6}}    \end{picture}}

\newcommand{\twthonth}{\begin{picture}(5,4)(0,-4)
\put(0,0){\line(1,0){5}}
\multiput(0,-1)(0,-1){3}{\line(1,0){2}}
\put(0,0){\line(0,-1){4}}
\multiput(1,0)(1,0){2}{\line(0,-1){3}}
\put(0,-4){\line(1,1){1}}  
\put(2,-3){\line(1,1){3}}    \end{picture}}

\newcommand{\ontwonon}{\begin{picture}(2,3)(0,-3)
\put(0,0){\line(1,0){2}}
\multiput(0,-1)(0,-1){2}{\line(1,0){1}}
\put(0,0){\line(0,-1){3}}
\put(1,0){\line(0,-1){2}}
\put(0,-3){\line(1,1){1}}
\put(1,-1){\line(1,1){1}}   \end{picture}}

\newcommand{\twtwtwtw}{\begin{picture}(4,4)(0,-4)
\put(0,0){\line(1,0){4}}
\multiput(0,-1)(0,-1){2}{\line(1,0){2}}
\put(0,0){\line(0,-1){4}}
\multiput(1,0)(1,0){2}{\line(0,-1){2}}
\put(0,-4){\line(1,1){4}}   \end{picture}}

\newcommand{\onthontw}{\begin{picture}(3,4)(0,-4)
\put(0,0){\line(1,0){3}}
\multiput(0,-1)(0,-1){3}{\line(1,0){1}}
\put(0,0){\line(0,-1){4}}
\put(1,0){\line(0,-1){3}}
\put(0,-4){\line(1,1){1}}
\put(1,-2){\line(1,1){2}} \end{picture}}

\newcommand{\twfoonfo}{\begin{picture}(6,5)(0,-5)
\put(0,0){\line(1,0){6}}
\multiput(0,-1)(0,-1){4}{\line(1,0){2}}
\put(0,0){\line(0,-1){5}}
\multiput(1,0)(1,0){2}{\line(0,-1){4}}
\put(0,-5){\line(1,1){1}}
\put(2,-4){\line(1,1){4}} \end{picture}}

\begin{document}

\title[Multipartitions, Generalized Durfee Squares and ...]
{Multipartitions, Generalized Durfee Squares and Affine Lie Algebra 
Characters}
\author{Peter Bouwknegt}
\address{Department of Physics and Mathematical Physics,
         University of Adelaide, Adelaide SA 5005, Australia}
\email{pbouwkne@physics.adelaide.edu.au}

\thanks{The author was supported by a 
 QEI$\!$I research fellowship from the Australian Research Council}
\dedicatory{Dedicated to Rodney J.~Baxter on his 60th birthday}

\begin{abstract}
We give some higher dimensional analogues of the Durfee square formula
and point out their relation to dissections of multipartitions.  We
apply the results to write certain affine Lie algebra characters in
terms of Universal Chiral Partition Functions.
\end{abstract}

\maketitle

\section{Introduction and background}

In this paper we will consider certain generalizations of an identity,
due to Euler, known as the Durfee square identity 
(see \cite{Andb} for an excellent introduction and historical account)
\begin{equation} \label{eqDURaa}
\frac{1}{\qn{\infty}} \eql \sum_{m\geq0} \frac{q^{m^2}}{\qn{m}\qn{m}}\,,
\end{equation}
where 
\begin{equation}  \label{eqDURab}
(z;q)_M \eql \prod_{k=1}^{M} \ (1-zq^{k-1}) \,,\qquad
\qn{M} \equiv (q;q)_M \,.
\end{equation}

There are various ways to prove this identity.  For instance, it follows 
as a limiting case of the $q$-analogue of Gauss' formula for the basic
hypergeometric series $_2\phi_1$ (see, e.g., \cite{Andc}).  
The most lucid proof, however,
employs the connection of \eqref{eqDURaa} to partitions \cite{Syl}
(see also \cite{HW,Andc}).  Henceforth we identify partitions 
$\bla = (\la_1,\la_2,\ldots)$, $\la_1\geq\la_2\geq \ldots
\geq 1$, and their graphical presentation in terms of Young diagrams 
\cite{Andc}
(see, e.g., Fig.\ 1.1 for the partition $\bla=(6,4,4,2)$).

Now, recall that
\begin{equation} \label{eqDURac}
(zq)_M^{-1} \eql \sum_{m,n\geq0}\ p_M(m,n)\ z^m q^n \,,
\end{equation}
where $p_{M}(m,n)$ denotes the number of partitions of $n$ into
$m$ parts in which no part exceeds $M$.  
In terms of Young diagrams, $p_M(m,n)$ is the number of diagrams
with $n$ boxes such that there are $m$ rows and no more than 
$M$ columns.


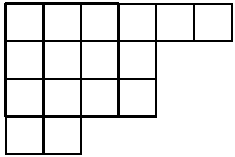
\begin{figure}
\begin{picture}(6,4)(0,0)
\put(0,0){\examtw}
\end{picture}
\caption{The partition $\bla=(6,4,4,2)$ and its $3\times3$ Durfee 
square}
\end{figure}


Thus, the left hand side of \eqref{eqDURaa} is clearly the generating 
function for all partitions, while each summand on the right hand 
side correspond to all partitions which fit at most an $m\times m$ 
`Durfee square' in the upper left hand corner of the Young diagram.
(The $3\times3$ Durfee square for the partition $\bla=(6,4,4,2)$
is indicated in Fig.\ 1.1.)
Summing over all $m$ clearly generates the total set of partitions
as well.  This proves \eqref{eqDURaa}.
In fact, by keeping track of the number of columns and rows 
in the above argument we have the following generalization of 
\eqref{eqDURaa} due to Cauchy
\begin{equation} \label{eqDURad}
\frac{1}{(zq)_M} \eql \sum_{m\geq0} \  \frac{q^{m^2}z^m}{\qny{z}{m}}
\qbin{M}{m} \,,
\end{equation}
where 
\begin{equation} \label{eqDURae}
\qbin{m}{n} \eql \frac{\qn{m}}{\qn{n} \qn{m-n}} \,,
\end{equation}
for $0\leq n\leq m$ (and zero otherwise), denotes the $q$-binomial 
(Gaussian polynomial).

Instead of dissecting partitions according to their maximal Durfee
square, Andrews considered dissections by (maximal) rectangles 
whose base to height ratio is $r:s$ and obtained the following
generalization of \eqref{eqDURad} \cite{Anda} 
\begin{equation} \label{eqDURaf} 
\frac{1}{(zq)_M} \eql 
\sum_{i,j}\ \sum_{m\geq0} \frac{q^{(rm+i)(sm+j)} z^{rm+i}}{
\qny{z}{sm+j-1+\delta_{i,0}+\delta_{i,r}} } 
\qbin{M+rm+i\delta_{j,s}-sm-j}{rm+i\delta_{j,s}} \,. 
\end{equation}
where the sum over $(i,j)$ is over all pairs
\begin{equation} \label{eqDURag}
(i,j) ~\in~ \{ (i,j)=(0,0) \ \text{or} \ 1\leq i\leq r,\ 
1\leq j \leq s,\ (i,j)\neq(r,s)\} \,.
\end{equation}
In fact, the identity (\ref{eqDURaf}) is valid even if $r$ and $s$ are
not relatively prime, as is obvious from Andrews' proof.
For $(r,s)=(1,1)$, Eq.\ \eqref{eqDURaf} reduces to \eqref{eqDURad}, 
while for $(r,s)=(2,1)$ it gives an identity which appears explicitly in
Ramanujan's lost notebook (see \cite{Andb}).

In this paper we will consider further generalizations of \eqref{eqDURaf} 
by considering simultaneous dissections of multipartitions. 
The resulting formulas are useful in deriving expressions for the 
chiral characters of 2D conformal field theories (in particular the
characters of modules of affine Lie algebras) in terms of so-called
universal chiral partition functions (UCPF's).

\section{Durfee systems}

We will be concerned with identities of the form
\begin{multline} \label{eqDURba} 
\frac{1}{\prod_i (z_iq)_{M_i}}  \eql  \sum_k\ 
\sideset{}{'}\sum_{\sumtop{\bm\in \ZZ_+^n}{\bn-\bK\cdot\bm=\bQ^{(k)}}}\  
\left( \prod_i z_i^{m_i+a^{(k)}_i} 
\right) q^{(\bm+\ba^{(k)})\cdot(\bn +\bb^{(k)})} 
  \frac{1}{\prod_i (z_iq)_{n_i}} \\
 \times \prod_i \qbin{ M_i+m_i-(n_i+b^{(k)}_i)}{m_i} \,.
\end{multline}
where $\bK\in GL(n,\QQ)$ is a symmetric matrix 
and the sum over $k$ is over a (finite) set of 
sectors.  
In each sector $k$, the sum over 
$\bm$ is over those $\bm\in(\ZZ_+)^n$ (here $\ZZ_+$ denotes the set 
of non-negative integers) such that $\bK\cdot\bm + 
\bQ^{(k)}\in (\ZZ_+)^n$, while
$\bn = \bK\cdot\bm + \bQ^{(k)}$. 

\begin{definition} A {\it Durfee system} 
for $\bK\in GL(n,\QQ)$, of length $L$, is a collection of
$n$-dimensional vectors, $(\bQ^{(k)},\ba^{(k)},\bb^{(k)})$,
$k=0,\ldots,L-1$, such that \eqref{eqDURba} is satisfied for all
$M_i\in \ZZ_+$ and $z_i$ ($i=1,\ldots,n$).
\end{definition}

Andrews' $(r,s)$-generalization of the classical Durfee formula,
discussed in Sect.\ 1, can now be formulated as

\begin{theorem} Let $r,s\in\NN$.  
A Durfee system of length $L=rs$, for the 
$1\times1$ matrix $\bK = s/r$, is given by
\begin{eqnarray} \label{eqDURbb} 
\bQ^{(i,j)} & \eql & j-1+\delta_{i,0} + \delta_{i,r} - \frac{s}{r} \,i
  \delta_{j,s} \,, \nonumber \\
\ba^{(i,j)} & \eql & i (1-\delta_{j,s}) \,, \nonumber \\
\bb^{(i,j)} & \eql & 1- \delta_{i,0} - \delta_{i,r} \,,
\end{eqnarray}
where $k=(i,j)$ runs over the $rs$ sectors as in \eqref{eqDURag}.
\end{theorem}

In the remainder of this paper 
we restrict ourselves to non-negative integer-valued, symmetric 
matrices $\bK$, i.e., $\bK\in GL(n,\ZZ_+)$, and Durfee systems
$(\bQ^{(k)},\ba^{(k)},\bb^{(k)})$ of $n$-vectors with entries 
in $\ZZ_+$.  In this case the sum in \eqref{eqDURba} is over all 
$m_i\geq0$ and $n_i\in\ZZ_+$ is determined by $\bn=\bK\cdot\bm +\bQ^{(k)}$.

Before giving examples, let us first explore some consequences of 
\eqref{eqDURba}.
By replacing $z_i\to z_i q^{p_i}$
in \eqref{eqDURba}, for some $\bp\in \ZZ^n$, 
using the expansion\footnote{Note that \eqref{eqDURbc}
itself can be interpreted as a length-1 Durfee system for the trivial
matrix $\bK=0$ with $(\bQ,\ba,\bb)=(0,0,1)$.}
\begin{equation} \label{eqDURbc}
\frac{1}{(zq)_M } \eql \sum_{m\geq0} \ (zq)^m \qbin{M+m-1}{m} \,,
\end{equation}
and shifting the summation variables, we find 
\begin{multline}\label{eqDURbd} 
\prod_i \qbin{M_i+N_i}{M_i} \eql  
\sum_k \sum_{\sumtop{\bm\in\ZZ_+^n}{ \bn - \bK\cdot \bm = \bQ^{(k)} + \bp}}
  \ q^{( \bm +\ba^{(k)})\cdot(\bn +\bb^{(k)}) } 
  \prod_i\qbin{M_i+m_i-(n_i+b^{(k)}_i)}{m_i} \\  
\times \qbin{N_i+n_i-(m_i+a^{(k)}_i)}{n_i} \,, 
\end{multline}
for arbitrary $\bp\in\ZZ^n$.  Note that in this formula the
summation variables $(\bm,\bn)$ appear on a more symmetrical
footing.

By taking the limit $M_i\to\infty$ in \eqref{eqDURba} we find
\begin{equation} \label{eqDURbe}
\frac{1}{\prod_i \qny{z_i}{\infty}} \eql 
\sum_{k} \sum_{\sumtop{\bm\in\ZZ_+^n}{  \bn - \bK\cdot \bm = \bQ^{(k)} } } 
  \ \left( \prod_i z_i^{m_i + a^{(k)}_i} \right) 
\frac{ q^{( \bm +\ba^{(k)})\cdot(\bn+\bb^{(k)})}}{
\prod_i \qn{m_i}\qny{z_i}{n_i} } \,,
\end{equation}
while by specializing \eqref{eqDURbe} to $z_i = q^{p_i}$, we find 
a generalization of the classical Durfee formula \eqref{eqDURaa}
\begin{equation} \label{eqDURbf}
\frac{1}{\qn{\infty}^n} \eql 
\sum_{k} \sum_{ \bn - \bK\cdot \bm = \bQ^{(k)} + \bp}  \ 
\frac{ q^{( \bm +\ba^{(k)})\cdot(\bn +\bb^{(k)})}}{\prod\qn{m_i}\qn{n_i} } 
\,,
\end{equation}
for any constant vector $\bp\in \ZZ^n$.  Of course, this equation can
also be obtained from \eqref{eqDURbd} by letting all $M_i\to\infty$.
Other interesting formulas are obtained by taking different specializations
of \eqref{eqDURbd}.  

The search for identities of the type \eqref{eqDURba} in 
dimension $n$ is greatly facilitated by using results in lower dimensions.
Indeed, by putting $z_i=0$ for some $i=i_0$ in \eqref{eqDURba}, the right 
hand side only receives contributions from the sectors $k$ for which 
$a_{i_0}^{(k)}=0$.  For those sectors only the term $m_{i_0}=0$ 
contributes in the summation, and \eqref{eqDURba} reduces to a similar
identity in dimension $n-1$.  Summarizing, if we know identities
for a $(n-1)\times (n-1)$ subblock of $\bK$, then we learn about 
the components $(Q^{(k)}_i,a^{(k)}_i,b^{(k)}_i)$, $i\neq i_0$, for all
sectors $k$ for which $a^{(k)}_{i_0}=0$.   

We now discuss the correspondence of Durfee systems with
multipartitions.  Suppose we have a Durfee system
$(\bQ^{(k)},\ba^{(k)},\bb^{(k)})$ for $\bK\in GL(n,\ZZ_+)$.
Consider Eq.\ \eqref{eqDURbf} for $\bp=0$.  
The left hand side is the generating 
series for all multipartitions $(\bla^{(1)},\bla^{(2)},
\ldots,\bla^{(n)})$.
Each term in the summand on the right hand side of \eqref{eqDURbf}
is a product (over $i$) of terms of the form 
\begin{equation} \label{eqDURbg}
\frac{q^{(m+a)(n+b)}}{\qn{m}\qn{n}} \,.
\end{equation}
By associating to \eqref{eqDURbg} a set of partitions of the form
indicated in Fig.\ 2.1, each term in the summand on the right hand 
side of \eqref{eqDURbf} is in 1--1 correspondence with a set of 
multipartitions.  


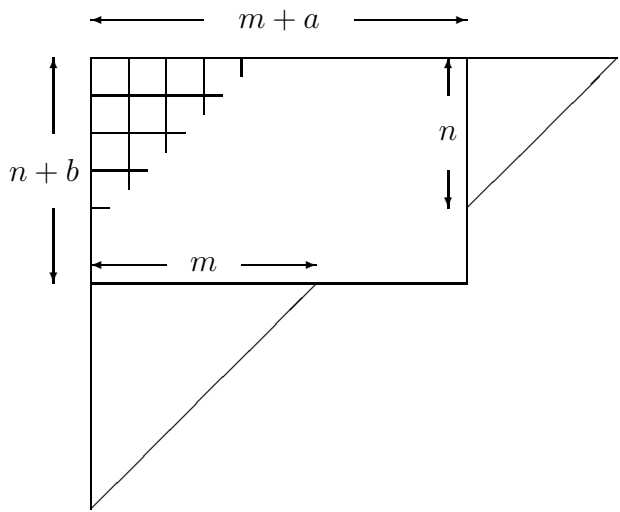
\begin{figure} 
\begin{picture}(24,14)(0,0)
\put(0,0){\mnab}
\end{picture}
\caption{Set of partitions with generating function \eqref{eqDURbg}}
\end{figure}


One possible strategy for proving the existence of a Durfee system is
therefore to show that the set of $n$-dimensional multipartitions
corresponding to the right hand side of \eqref{eqDURbf} is
non-overlapping and exhaustive.  By keeping track of the number of
rows and columns in each partition $\bla^{(i)}$, the generalization
\eqref{eqDURba} then easily follows.

After discussing some examples of Durfee systems in the following
sections we will explore some further consequences in the context of
affine Lie algebra characters.

\section{Examples}

In this section we will consider some examples of Durfee systems.

\begin{theorem} Consider the matrix $\bK\in GL(2,\ZZ_+)$ given by 
\begin{equation} \label{eqDURcaa}
\bK \eql \left( \matrix 1 & 1 \\ 1 & 2 \endmatrix \right) \,.
\end{equation}
We have a Durfee system $(\bQ^{(k)},\ba^{(k)},\bb^{(k)})$ for $\bK$
given by 
\begin{align}
\bQ^{(0)} & \eql \begin{pmatrix} 0 \\  0 \end{pmatrix}\,, &
\ba^{(0)} & \eql \begin{pmatrix} 0 \\  0 \end{pmatrix}\,, &
\bb^{(0)} & \eql \begin{pmatrix} 0 \\  0 \end{pmatrix}\,, \nonumber \\
\bQ^{(1)} & \eql \begin{pmatrix} 0 \\  1 \end{pmatrix}\,, &
\ba^{(1)} & \eql \begin{pmatrix} 0 \\  1 \end{pmatrix}\,, &
\bb^{(1)} & \eql \begin{pmatrix} 1 \\  0 \end{pmatrix}\,.
\end{align}
\end{theorem}

Let us illustrate, in some detail, how one might arrive at this result.
The $k=0$ term in \eqref{eqDURbf} (for $\bp=0$) is explicitly given by
\begin{equation} \label{eqDURla}
\sum_{\sumtop{n_1 -(m_1+m_2) = 0}{n_2 -(m_1+2m_2) = 0} }
\frac{q^{n_1m_1+n_2m_2}}{\qn{n_1}\qn{n_2}\qn{m_1}\qn{m_2}} \,.
\end{equation}
The set of bipartitions $(\bla^{(1)},\bla^{(2)})$
associated to \eqref{eqDURla}, according to the prescription of
Sect.\ 2, is depicted in Fig.\ 3.1 for low values of $\bm=(m_1,m_2)$.


\begin{figure}
\begin{picture}(30,7)(-0.5,-3.5)

\put(0,0){\oval(1,2)[l]}
\put(-0.5,-0.5){\zezezeze}
\put(1.5,-0.5){\zezezeze}
\put(2,0){\oval(1,2)[r]}
\put(1,-3){\makebox(0,0){$\bm=$(0,0)}}

\put(9,0){\oval(1,3)[l]}
\put(9,-1){\onononon}
\put(13,0){\zeonzeon}
\put(14.5,0){\oval(1,3)[r]}
\put(11.75,-3){\makebox(0,0){(1,0)}}

\put(20,0){\oval(1,4)[l]}
\put(20,0.5){\zeonzeon}
\put(23.5,-1.5){\ontwontw}
\put(26.5,0){\oval(1,4)[r]}
\put(23.25,-3){\makebox(0,0){(0,1)}}

\end{picture}
\begin{picture}(30,8)(-0.5,-3.5)

\put(0,0){\oval(1,5)[l]}
\put(0,-1){\ontwontw}
\put(5,-2){\onthonth}
\put(9,0){\oval(1,5)[r]}
\put(4.5,-3.5){\makebox(0,0){(1,1)}}

\put(12,0){\oval(1,5)[l]}
\put(12,-2){\twtwtwtw}
\put(18,0){\zetwzetw}
\put(20.5,0){\oval(1,5)[r]}
\put(16.25,-3.5){\makebox(0,0){(2,0)}}

\put(25,0){\makebox(0,0){$\ldots$}}

\end{picture}
\caption{The $k=0$ sector}
\end{figure}


Clearly these do not exhaust the set of all bipartitions.  For instance,
if $\bla^{(1)}=\emptyset$ (indicated by a $\bullet$ in Fig.\ 3.1)
and $\bla^{(2)}\neq\emptyset$, then $\bla^{(2)}$ necessarily has two 
or more rows.  Thus, the set of bipartitions depicted in 
Fig.\ 3.2 is missing in \eqref{eqDURla}.


\begin{figure}
\begin{picture}(4,2)(0,-1)

\put(0,0){\oval(1,2)[l]}
\put(-0.5,-0.5){\zezezeze}
\put(1.5,-0.5){\ononzeon}
\put(3.5,0){\oval(1,2)[r]}

\end{picture}
\caption{Some missing bipartitions}
\end{figure}
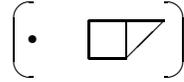


If this set of bipartitions is to be included as the $\bm=(0,0)$
term of another sector, say $k=1$, then this immediately fixes 
all components of $(\bQ^{(1)},\ba^{(1)},\bb^{(1)})$ with the exception
of $b^{(1)}_1$. [Note that this component is also unconstrained by
consideration of the two $1\times1$ subblocks of $\bK$, as discussed 
in Sect.\ 2.]
Consideration of the $\bm=(1,0)$ term in the $k=1$ sector, however,
uniquely fixes $b^{(1)}_1$ as well and we arrive at the conclusion that 
\eqref{eqDURla} needs to be supplemented by 
\begin{equation} \label{eqDURlb}
\sum_{\sumtop{n_1 -(m_1+m_2) = 0}{n_2 -(m_1+2m_2) = 1} }
\frac{q^{(n_1+1)m_1+n_2(m_2+1)}}{\qn{n_1}\qn{n_2}\qn{m_1}\qn{m_2}} \,. 
\end{equation} 
The set of bipartitions
in the $k=1$ sector, arising from \eqref{eqDURlb}
for low values of $\bm$, is depicted in Fig.\ 3.3.


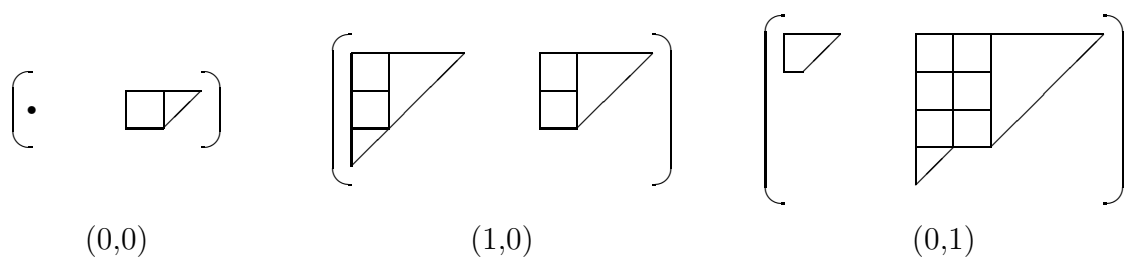
\begin{figure}
\begin{picture}(30,7)(-0.5,-3.5)

\put(0,0){\oval(1,2)[l]}
\put(-0.5,-0.5){\zezezeze}
\put(1.5,-0.5){\ononzeon}
\put(3.5,0){\oval(1,2)[r]}
\put(1.75,-3){\makebox(0,0){(0,0)}}

\put(9,0){\oval(1,4)[l]}
\put(9,-1.5){\ontwonon}
\put(13,-0.5){\ontwzetw}
\put(16,0){\oval(1,4)[r]}
\put(12.5,-3){\makebox(0,0){(1,0)}}

\put(20,0){\oval(1,5)[l]}
\put(20,1){\zeonzeon}
\put(23.5,-2){\twthonth}
\put(28.5,0){\oval(1,5)[r]}
\put(24.25,-3){\makebox(0,0){(0,1)}}

\end{picture}

\begin{picture}(30,8)(-0.5,-4)

\put(0,0){\oval(1,6)[l]}
\put(0,-1.5){\onthontw}
\put(5,-2.5){\twfoonfo}
\put(11,0){\oval(1,6)[r]}
\put(5.5,-4){\makebox(0,0){(1,1)}}

\put(15,0){\makebox(0,0){$\ldots$}}

\end{picture}
\caption{The $k=1$ sector}
\end{figure}


Together, the sets of bipartitions of Figs.\ 3.1 and 3.3 are seen 
to be non-overlapping and to exhaust the set of all bipartitions,
at least to low order, so it seems that no other sectors are required.
The proof that this works to all orders requires a bit more work and
will be omitted. 

A slightly more complicated Durfee system is given in 
\begin{theorem} Let
\begin{equation} \label{eqDURca}
\bK \eql \left( \matrix 2 & 1 \\ 1 & 2 \endmatrix \right) \,.
\end{equation}
The following constitutes a Durfee system for $\bK$
\begin{align} \label{eqDURcb}
\bQ^{(0)} & \eql \begin{pmatrix} 0 \\ 0 \end{pmatrix} \,, & 
\ba^{(0)} & \eql \begin{pmatrix} 0 \\ 0 \end{pmatrix} \,, & 
\bb^{(0)} & \eql \begin{pmatrix} 0 \\ 0 \end{pmatrix} \,, \nonumber \\
\bQ^{(1)} & \eql \begin{pmatrix} 0 \\ 1 \end{pmatrix} \,, & 
\ba^{(1)} & \eql \begin{pmatrix} 0 \\ 1 \end{pmatrix} \,, & 
\bb^{(1)} & \eql \begin{pmatrix} 0 \\ 0 \end{pmatrix} \,, \nonumber \\
\bQ^{(2)} & \eql \begin{pmatrix} 1 \\ 1 \end{pmatrix} \,, & 
\ba^{(2)} & \eql \begin{pmatrix} 1 \\ 0 \end{pmatrix} \,, & 
\bb^{(2)} & \eql \begin{pmatrix} 0 \\ 0 \end{pmatrix} \,. 
\end{align}
\end{theorem}

The reasoning parallels that of Theorem 3.1.  The first few sets of
contributing bipartitions, for the sectors $k=0,1,2$, are depicted in 
Figs.\ 3.4--3.6, respectively. 


\begin{figure}
\begin{picture}(30,6)(-0.5,-3)

\put(0,0){\oval(1,2)[l]}
\put(-0.5,-0.5){\zezezeze}
\put(1.5,-0.5){\zezezeze}
\put(2,0){\oval(1,2)[r]}
\put(1,-3){\makebox(0,0){(0,0)}}

\put(9,0){\oval(1,4)[l]}
\put(9,-1.5){\ontwontw}
\put(14,0.5){\zeonzeon}
\put(15.5,0){\oval(1,4)[r]}
\put(12.25,-3){\makebox(0,0){(1,0)}}

\put(20,0){\oval(1,4)[l]}
\put(20,0.5){\zeonzeon}
\put(23.5,-1.5){\ontwontw}
\put(26.5,0){\oval(1,4)[r]}
\put(23.25,-3){\makebox(0,0){(0,1)}}

\end{picture}

\begin{picture}(30,9)(-0.5,-4.5)

\put(0,0){\oval(1,5)[l]}
\put(0,-2){\onthonth}
\put(6,-2){\onthonth}
\put(10,0){\oval(1,5)[r]}
\put(5,-4.5){\makebox(0,0){(1,1)}}

\put(12,0){\oval(1,7)[l]}
\put(12,-3){\twfotwfo}
\put(20,1){\zetwzetw}
\put(22.5,0){\oval(1,7)[r]}
\put(17.25,-4.5){\makebox(0,0){(2,0)}}

\put(28,0){\makebox(0,0){$\ldots$}}
\end{picture}

\caption{The $k=0$ sector}
\end{figure}
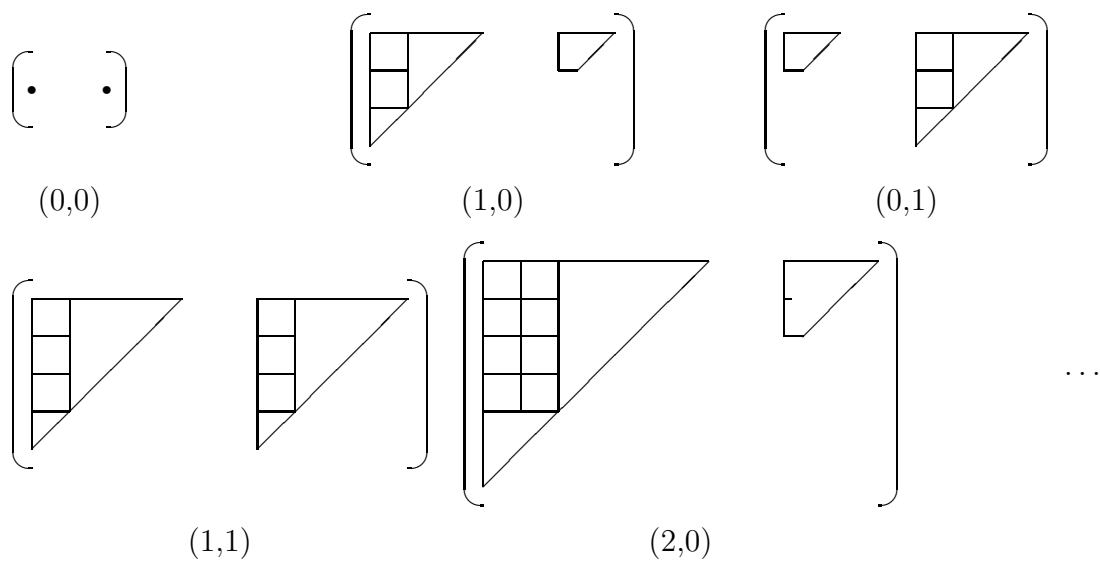


\begin{figure}
\begin{picture}(30,7)(-0.5,-3.5)

\put(0,0){\oval(1,2)[l]}
\put(-0.5,-0.5){\zezezeze}
\put(2.5,-0.5){\ononzeon}
\put(4.5,0){\oval(1,2)[r]}
\put(2.25,-3.5){\makebox(0,0){(0,0)}}

\put(8.5,0){\oval(1,4)[l]}
\put(8.5,-1.5){\ontwontw}
\put(13.5,-0.5){\ontwzetw}
\put(16.5,0){\oval(1,4)[r]}
\put(12.5,-3.5){\makebox(0,0){(1,0)}}

\put(20,0){\oval(1,5)[l]}
\put(20,1){\zeonzeon}
\put(23.5,-2){\twthonth}
\put(28.5,0){\oval(1,5)[r]}
\put(24.25,-3.5){\makebox(0,0){(0,1)}}

\end{picture}
\caption{The $k=1$ sector}
\end{figure}


\begin{figure}
\begin{picture}(30,7)(-0.5,-3.5)

\put(0,0){\oval(1,2)[l]}
\put(0,-0.5){\ononzeon}
\put(4,-0.5){\zeonzeon}
\put(5.5,0){\oval(1,2)[r]}
\put(2.25,-3.5){\makebox(0,0){(0,0)}}

\put(8.5,0){\oval(1,5)[l]}
\put(8.5,-2){\twthonth}
\put(15.5,0){\zetwzetw}
\put(18,0){\oval(1,5)[r]}
\put(13.25,-3.5){\makebox(0,0){(1,0)}}

\put(20,0){\oval(1,5)[l]}
\put(20,0){\ontwzetw}
\put(25,-2){\onthonth}
\put(29,0){\oval(1,5)[r]}
\put(24.5,-3.5){\makebox(0,0){(0,1)}}

\end{picture}
\caption{The $k=2$ sector}
\end{figure}
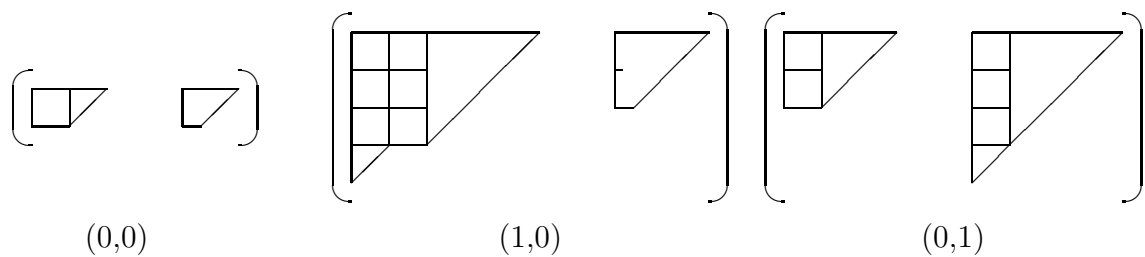


Theorem 3.2 has the following higher dimensional generalization
\begin{theorem} Let $\bK\in GL(n,\ZZ_+)$ be defined by
\begin{equation} \label{eqDURmf}
\bK \eql \begin{pmatrix} 2 & 1 & 1 & \ldots & 1 \\
                         1 & 2 & 1 & \ldots & 1 \\
                         \vdots & & \ddots & & \vdots \\
                         \vdots & &  & & \vdots \\
                         1 & 1 &  1 & \ldots & 2 \end{pmatrix} \,.
\end{equation}
We have a Durfee system of length $L=n+1$, given by the $n$-vectors
\begin{eqnarray} \label{eqDURga}
\bQ^{(k)} & \eql & (\underbrace{0,0,\ldots,0}_{n-k},
\underbrace{1,\ldots,1}_{k}) \,, \nonumber\\
\ba^{(k)} & \eql & (\underbrace{0,\ldots,0}_{n-k},1, 
\underbrace{0,\ldots,0}_{k-1}) \,, \nonumber\\
\bb^{(k)} & \eql & (\underbrace{0,0,\ldots,0}_n) \,, 
\end{eqnarray}
for $k=0,\ldots,n$. 
\end{theorem}

\begin{remark} Note that the length of the Durfee system in Theorem 3.3
is given by $L=n+1=\det\bK$.  We believe this is a general feature of
Durfee systems for which $\bb^{(k)}=0$ for all $k=0,\ldots,L-1$
(see also the discussion in Sect.\ 6).
\end{remark}

\section{Shift operation}

It turns out that, once a Durfee system for some $\bK\in GL(n,\ZZ_+)$
has been established, it is rather straightforward to obtain a 
Durfee system for a class of deformations of $\bK$.  These 
deformations are given in terms of 
a ``charge vector'' $\bt=(t_1,\ldots,t_n)$, $t_i\in\ZZ_+$, and
a positive integer $M\in\ZZ_+$ as\footnote{These deformations were motivated
by the ``shift operation'' on $K$-matrices describing fractional 
quantum Hall systems (see \cite{ABS} and references therein).}
\begin{equation} \label{eqDURha}
\bK_{M,\bt} \eql \bK + M\,\bt\bt^T\,.  
\end{equation}
For instance, consider the deformation $\bK_{M,\bt}$ of the two-dimensional 
identity matrix
\begin{equation} \label{eqDURhb}
\bK_{M,\bt} \eql \begin{pmatrix}
t_1^2 M +1 & t_1t_2M \\
t_1t_2M & t_2^2 M +1 \end{pmatrix} \,,
\end{equation}
where we can assume that $t_1\leq t_2$.  Note that 
the matrix $\bK$ of Eq.\ \eqref{eqDURca}
is of this form with $M=1$, $\bt=(1,1)$.

\begin{theorem} The matrix $\bK_{M,\bt}$ of Eq.\ \eqref{eqDURhb}
admits a length $L=(t_1^2 + t_2^2)M +1$ Durfee system.  There are 
$t_2^2M$ sectors given by 
\begin{multline}
\bQ \eql \underbrace{
         \left( \matrix t_1^2M+1 \\ t_2^2M \endmatrix \right),\ 
         \left( \matrix t_1^2M+1 \\ t_2^2M-1 \endmatrix \right),\ldots,
         \left( \matrix t_1^2M+1 \\ t_1^2M+2 \endmatrix \right)}_{
                     (t_2^2-t_1^2)M}    ,\\ 
         \underbrace{
         \left( \matrix t_1^2M \\ t_1^2M+1 \endmatrix \right),\ 
         \left( \matrix t_1^2M-1 \\ t_1^2M \endmatrix \right),\ldots,
         \left( \matrix  0 \\ 1 \endmatrix \right)}_{t_1^2M} \,,
\end{multline}
with
$\ba=\begin{pmatrix} 0 \\ 0 \end{pmatrix}$, 
$\bb=\begin{pmatrix} 0 \\ 0 \end{pmatrix}$,
$t_1^2M$ sectors
given by
\begin{equation}
\bQ \eql \left( \matrix t_1^2M \\ t_1^2M \endmatrix \right),\ 
         \left( \matrix t_1^2M-1 \\ t_1^2M-1 \endmatrix \right),\ldots,
         \left( \matrix  1 \\ 1 \endmatrix \right) \,,
\end{equation}
with 
$\ba=\begin{pmatrix} 1 \\ 0 \end{pmatrix}$, 
$\bb=\begin{pmatrix} 0 \\ 0 \end{pmatrix}$,
and and the `vacuum sector' 
$\bQ =\begin{pmatrix} 0 \\ 0 \end{pmatrix}$,
$\ba=\begin{pmatrix} 0 \\ 0 \end{pmatrix}$, 
$\bb=\begin{pmatrix} 0 \\ 0 \end{pmatrix}$.
\end{theorem}

For deformations \eqref{eqDURha}, with $\bK=\id$, we have 
\begin{equation} \label{eqDURhc}
\det \bK_{M,\bt}  \eql  (\bt^T\cdot \bt)\, M + 1\,,
\end{equation}
which can be written as 
\begin{equation} \label{eqDURhd}
\det \bK_{M,\bt}  \eql  \tr(\bK_{M,\bt}-\id) + 1\,.
\end{equation}
In fact, if $n=2$, the matrix $\bK_{M,\bt} = \id + M \, \bt \bt^T$
is the most general symmetric,
non-negative integer-valued matrix satisfying \eqref{eqDURhd}.
Note that the length of the Durfee system in Theorem 4.1 is
again given by $\det \bK_{M,\bt}$.

\section{The UCPF and character identities}

Consider the  ``Universal Chiral Partition Function'' (UCPF) 
(see \cite{BM} and references therein)
\begin{equation} \label{eqDURea}
Z(\bK;\bQ,\bu|\bz;q) 
\eql \sum_{\bm\in\ZZ_+^n} \left( \prod_i z_i^{m_i} \right)\ 
 q^{\frac{1}{2} \bm\cdot\bK\cdot\bm
 + \bQ\cdot \bm} \prod_i \qbin{ ((\id-\bK)\cdot \bm + \bu)_i}{m_i} \,,
\end{equation}
where $\bK\in GL(n,\ZZ_+)$, $Q_i\in\ZZ_+$ and $u_i\in\ZZ_+\cup\{\infty\}$, 
$i=1,\ldots,n$.\footnote{The considerations in this section
can easily be generalized to triples $(\bK;\bQ,\bu)$ with entries in
$\QQ$, provided appropriate restrictions on the summation variables
$m_i$ in \eqref{eqDURea} are made.}

The following theorem is derived by elementary algebra
\begin{theorem} Assume that 
$(\bQ^{(k)},\ba^{(k)},\bb^{(k)})$ forms a Durfee 
system for a symmetric $\bK\in GL(n,\ZZ_+)$.  Define 
\begin{equation}  \label{eqDURec}
\bQ^{\prime(k)} \eql - \bK^{-1}\cdot \bQ^{(k)} \,, \qquad
z'_i \eql \prod_j \, z_j^{-K_{ij}} \,.			
\end{equation}
Then we have the following identity
\begin{multline} \label{eqDURka}
\sum_k   \left( \prod_i z_i^{-Q^{(k)}_i} \right)
 q^{\frac{1}{2} \bQ^{(k)} \cdot\bK^{-1}\cdot\bQ^{(k)} + 
  \ba^{(k)}\cdot \bb^{(k)}} 
 Z(\bK;\bQ^{(k)} + \bb^{(k)},\bM - (\bQ^{(k)} + \bb^{(k)})|\bz';q)\\   
 \times 
 Z(\bK^{-1}; \bQ^{\prime(k)} + \ba^{(k)},\bN - (\bQ^{\prime(k)} + \ba^{(k)})|
 \bz;q) \\ \eql
\sum_{\bp\in\ZZ^n} \left( \prod_i z_i^{p_i} \right) \ 
 q^{\frac{1}{2} \bp\cdot \bK^{-1}\cdot \bp} 
 \prod_i \qbin{M_i + N_i + ((\id-\bK^{-1})\cdot \bp)_i}{M_i+p_i} \,.
\end{multline}
for all $\mathbf{M}, \mathbf{N} \in \ZZ_+^n$.
\end{theorem}

\begin{remark}
Note that the polynomials 
$P^{(k)}_{\bM}(\bz;q) \equiv 
Z(\bK;\bQ^{(k)} + \bb^{(k)},\bM - (\bQ^{(k)} + \bb^{(k)})|\bz';q)$
and $Q^{(k)}_{\bN}(\bz;q) \equiv 
Z(\bK^{-1};\bQ^{'(k)} + \ba^{(k)},\bN - (\bQ^{'(k)} + \ba^{(k)})|\bz;q)$,
entering Eq.\ \eqref{eqDURka}, 
all arise as a solution to the same (i.e.\ $k$-independent) set of 
recursion relations ($i=1,\ldots,n$) \cite{ABS}
\begin{eqnarray} \label{eqDURed}
P_{\bM}(\bz';q) & \eql & P_{\bM - \bfe_i}(\bz';q) + 
  z_i' q^{-\frac{1}{2} K_{ii} + M_i} P_{\bM-\bK\cdot \bfe_i}(\bz';q)  \,,
 \nonumber \\
Q_{\bN}(\bz;q) & \eql & Q_{\bN - \bfe_i}(\bz;q) + 
  z_i q^{-\frac{1}{2} K^{-1}_{ii} + N_i} 
  Q_{\bN-\bK^{-1}\cdot \bfe_i}(\bz;q)  \,,
\end{eqnarray}
where $\bfe_i$ denotes the unit vector in the $i$-direction and where
we have used
\begin{equation}
\qbin{M}{m} \eql \qbin{M-1}{m} + q^{M-m} \qbin{M-1}{m-1} \,.
\end{equation}
\end{remark} \medskip

For the application of Theorem 5.1 to affine Lie algebra characters
let us consider the limiting form of \eqref{eqDURea} as 
$\bu\to\infty$, i.e.,
\begin{equation}
Z_\infty(\bK;\bQ|\bz,q) \eql \lim_{{\bu}\to\infty} 
Z(\bK;\bQ,\bu|\bz,q) \eql
\sum_{\bm} \left( \prod_i z_i^{m_i} \right) 
\frac{q^{\frac{1}{2} \bm\cdot\bK\cdot\bm + \bQ\cdot \bm}}{ \prod_i
 \qn{m_i} } \,.
\end{equation}

\begin{remark} The limiting UCPF's are not all independent.  For instance,
by using the simple relation
\begin{equation}
\frac{1}{\qn{m}} - \frac{q^m}{\qn{m}} \eql \frac{1}{\qn{m-1}} \,.
\end{equation}
we find 
\begin{equation} \label{eqDUReaf}
Z_\infty(\bK;\bQ) \eql Z_\infty(\bK;\bQ+{\bfe}_i) + z_i q^{\frac{1}{2}
 {\bfe}_i \cdot \bK \cdot {\bfe}_i +  \bQ\cdot {\bfe}_i}
 Z_\infty(\bK;\bQ+\bK\cdot {\bfe}_i) \,.
\end{equation}
\end{remark} \medskip

By taking $\mathbf M\to\infty$ in \eqref{eqDURka} we obtain
\begin{corollary} Let $(\bQ^{(k)},\ba^{(k)},\bb^{(k)})$ be a Durfee
system for $\bK\in GL(n,\ZZ_+)$ of length $L$. Define 
$\bQ^{'(k)}$ and $z_i'$ by Eq.\ \eqref{eqDURec}.  We then have
\begin{multline} \label{eqDURma}
\sum_{k=0}^{L-1}  \left( \prod_i z_i^{-Q^{(k)}_i} \right)
 q^{\frac{1}{2} \bQ^{(k)} \cdot\bK^{-1}\cdot\bQ^{(k)} + 
  \ba^{(k)}\cdot \bb^{(k)}} 
Z_\infty(\bK;\bQ^{(k)} + \bb^{(k)}|\bz';q) \\ \times
Z_\infty(\bK^{-1}; \bQ^{\prime(k)} + \ba^{(k)}|\bz;q) 
 \eql \frac{1}{\qn{\infty}^n}
 \sum_{\bp\in\ZZ^n} \left( \prod_i z_i^{p_i} \right) \ 
 q^{\frac{1}{2} \bp\cdot \bK^{-1}\cdot \bp} \,.
\end{multline}
\end{corollary}

Now suppose that the bilinear form $\bp\cdot\bK^{-1}\cdot\bp$ is
chosen in such a way that it equals the standard bilinear form on the
weight lattice $\Lambda_{\text{w}}$
of a simple Lie algebra $\fg$ of rank $n$ and that
the sum over $\bp\in \ZZ^n$ corresponds to the sum over the weight
lattice.  Then, provided $\fg$ is simply-laced, the right hand
side of \eqref{eqDURma} can be recognized as the Frenkel-Kac character
of the sum of the level-1 integrable highest weight modules of the
affine Lie algebra $\wh{\fg}$ (see, e.g., \cite{Kac})%
\footnote{The irreducible characters can be recovered by suitably 
restricting the sum over $\bp$.}
Thus, in such cases, Corollary 5.2 provides an expression for the 
level-1 characters of $\wh{\fg}$ in terms of UCPF's based
on the bilinear form constructed out of $\bK\oplus\bK^{-1}$.
This has important applications in the study of  
quasiparticles in the conformal field 
theory descriptions of certain non-Abelian fractional quantum Hall
states \cite{ABGS,ABS}.  In fact, 
these applications were the main motivation
for the present study.

As an example, consider $\fg=\sln$.  The weights $\{\bep_1,\ldots,
\bep_{n+1}\}$,
of the fundamental $(n+1)$-dimensional 
representation $L(\Lambda_1)$ of $\sln$ satisfy 
$\bep_i\cdot \bep_j = \delta_{ij} - 1/(n+1)$.  A suitable basis
of the weight lattice
$\Lambda_{\text{w}}$ is given by the $\bep_i$, $i=1,\ldots,n$
(see Fig.\ 5.1 for $\mathfrak{sl}_3$).
Now note that  
\begin{equation}
(\sum_i p_i \bep_i) \cdot (\sum_j p_j \bep_j) 
\eql \bp \cdot \bK^{-1} \cdot \bp \,,
\end{equation}
where $\bK^{-1}$ is given by
\begin{equation}
\bK^{-1} \eql \frac{1}{n+1}
          \begin{pmatrix} n & -1 & -1 & \ldots & -1 \\
                         -1 & n & -1 & \ldots & -1 \\
                          \vdots & & \ddots & & \vdots \\
                           \vdots & &  & & \vdots \\
                         -1 & -1 &  -1 & \ldots & n \end{pmatrix} \,,
\end{equation}
which has an inverse $\bK$ given by Eq.\ \eqref{eqDURmf}.  
The ``dual sector'',
defined by $\bK$, corresponds to a particular basis of the root lattice 
of $\sln$ (see Fig.\ 5.1 for $\mathfrak{sl}_3$).
The weights of this basis are determined by \eqref{eqDURec}.


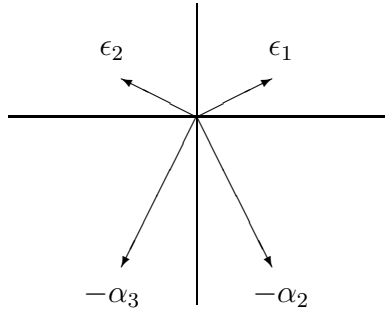
\begin{figure}
\begin{picture}(10,8)(-5,-5)
\put(-5,0){\line(1,0){10}}
\put(0,-5){\line(0,1){8}}
\put(0,0){\vector(2,1){2}}
\put(0,0){\vector(-2,1){2}}
\put(0,0){\vector(1,-2){2}}
\put(0,0){\vector(-1,-2){2}}
\put(2.25,-4.75){\makebox(0,0){$-\alpha_2$}}
\put(-2.25,-4.75){\makebox(0,0){$-\alpha_3$}}
\put(2.25,1.75){\makebox(0,0){$\epsilon_1$}}
\put(-2.25,1.75){\makebox(0,0){$\epsilon_2$}}
\end{picture}
\caption{$\mathfrak{sl}_3$ weights and roots}
\end{figure}


Thus, the sum over $\bp\in\ZZ^{n}$ is precisely over the weight
lattice of $\sln$ and combining Theorem 3.3 and Corollary 5.2 gives us
an expression for the character of the (sum over all) level-1
integrable highest weight modules of $\whsln$.  As a consistency
check, note that
\begin{equation}
\frac{1}{2} \bQ^{(k)} \cdot \bK^{-1} \cdot \bQ^{(k)} 
\eql \frac{k(n+1-k)}{2(n+1)} \,,\qquad k=0,\ldots,n\,,
\end{equation}
is indeed precisely the conformal dimension of the level-1 integrable 
highest weight module $L(\Lambda_k)$ of $\whsln$.

\section{Discussion and conclusions}

In this paper we have introduced higher dimensional analogues of the
classical Durfee square formula \eqref{eqDURaa}
in the form of ``Durfee systems'', we
explained their correspondence to multipartitions, and gave a few
examples.  We have also remarked on the application of Durfee systems,
in particular with regards to writing (chiral) characters of
two-dimensional conformal field theories in UCPF form.

A number of obvious questions come to mind.  Firstly, for which
symmetric $\bK \in GL(n,\ZZ_+)$ is it possible to find a Durfee
system?  It seems that this class of matrices is quite big.  In fact,
examples suggest that, provided $\det\bK\geq0$, a Durfee system always
exists (see \eqref{eqDURbc} for an example with $\det\bK=0$).
Secondly, how unique are Durfee systems for a given matrix $\bK$?
Clearly they are not unique. For instance, in the case of $\bK=s/r$
(see Theorem 2.2) we can construct Durfee systems of length $L=m^2 rs$
for all $m\in\NN$ by taking $(r,s)\to (mr,ms)$ in Eqs.\
\eqref{eqDURaf} and \eqref{eqDURag}.  Similar constructions exist for
the higher dimensional cases.  Another source of non-uniqueness
originates from possible symmetries of the matrix $\bK$.  For example,
interchanging the components of all vectors
$(\bQ^{(k)},\ba^{(k)},\bb^{(k)})$ in Theorem 3.2, provides another
Durfee system due to the $\ZZ_2$ permutation symmetry of the matrix
$\bK$ in \eqref{eqDURca}.
Thirdly, for a given $\bK$, what is the minimal length $L_{\text{min}}$
of a Durfee
system?  It seems that a special role is played by matrices for
which $L_{\text{min}}=\det\bK$, which seem to be closely related 
to matrices for which it is possible to choose a Durfee system 
for which $\bb^{(k)}=0$ for all $k$.  A large class of such matrices 
is provided by the shift 
deformations $\bK_{M,\bt}$ 
of the identity (see Eq.\ \eqref{eqDURha}) and, 
at least in two dimensions, it appears that
such deformations exhaust all matrices $\bK$ for which 
$L_{\text{min}}=\det\bK$.
Finally, is it possible to give a more `geometric' construction of
the vectors $(\bQ^{(k)},\ba^{(k)},\bb^{(k)})$?  Again,
in the case of matrices $\bK$ for which $L_{\text{min}}=\det\bK$ it seems 
that the set of $\bQ^{(k)}$ is given by a set of coset representatives
(with minimal non-negative components) of $\ZZ^n$ modulo the 
equivalences $\bm \sim \bm + \bK\cdot\bfe_i$ ($i=1,\ldots,n$).  
Note that in the case 
of \eqref{eqDURmf} the equivalence preserves the $\ZZ_{n+1}$ charge
$q=\sum im_i\ (\text{mod}\ n+1)$ of $\bm$ (``$n$-ality'') and that we find 
one coset representative for each $q\in\ZZ_{n+1}$.



\begin{thebibliography}{X}

\bibitem[1]{Anda}
G.E.~Andrews, {\it Generalizations of the Durfee square},
J.\ London Math.\ Soc.\ {\bf 3} (1971) 563-570.

\bibitem[2]{Andb}
G.E.~Andrews, {\it Partitions: Yesterday and today}, 
(New Zealand Mathematical Society, Wellington, 1979).

\bibitem[3]{Andc}
G.E.~Andrews, {\it The theory of partitions},
Encycl.\ of Math.\ and its Appl., Vol {\bf 2},
(Addison-Wesley, Reading, 1976).

\bibitem[4]{ABGS}
E.~Ardonne, P.~Bouwknegt, S.~Guruswamy and K.~Schoutens, 
{\it $K$-matrices for non-Abelian quantum Hall states},
Phys.\ Rev.\ {\bf B}, to appear, [{\tt cond-mat/9908285}].

\bibitem[5]{ABS}
E.~Ardonne, P.~Bouwknegt and K.~Schoutens, {\it Non-Abelian 
quantum Hall states -- exclusion statistics, $K$-matrices and duality --},
in preparation.

\bibitem[6]{BM}
A.~Berkovich and B.~McCoy, 
{\it The universal chiral partition function for exclusion statistics},
in ``Statistical Physics on the Eve of the 21st Century'',
Series on Adv.\ in Stat.\ Mech., Vol.\ {\bf 14}, pp 240-256,
eds.\ M.T.~Batchelor and L.T.~Wille,
(World Scientific, Singapore, 1999),
[{\tt hep-th/9808013}].

\bibitem[7]{HW}
G.H.~Hardy and E.M.~Wright, {\it An introduction to the 
theory of numbers}, (Oxford University Press, Oxford, 1960).

\bibitem[8]{Kac}
V.G.~Kac, {\it Infinite dimensional Lie algebras}, 
(Cambridge University Press, Cambridge, 1985).

\bibitem[9]{Syl}
J.J.~Sylvester, {\it A constructive theory of partitions, arranged in
three acts, an interact and an exodion}, Collected works, Vol.~{\bf 4}
(Cambridge University Press, Cambridge, 1912).


\end{thebibliography}
\end{document}